\documentclass[11.05pt,a4paper]{amsart}
\usepackage{amssymb,amsmath}
\usepackage{latexsym}
\usepackage{amsthm,amsfonts,amssymb,mathrsfs}
\usepackage{rotating}
\usepackage[leqno]{amsmath}
\usepackage{xspace}
\usepackage[all]{xy}
\usepackage{longtable}

\headsep=1cm \numberwithin{equation}{section}
\newtheorem{theorem}{Theorem}[section]
\newtheorem{lemma}[theorem]{Lemma}

\newtheorem{corollary}[theorem]{Corollary}

\theoremstyle{definition}
\newtheorem{definition}[theorem]{Definition}
\theoremstyle{remark}
\newtheorem{remark}[theorem]{Remark}
\newtheorem{example}[theorem]{Example}

\newcommand{\Ass}{\operatorname{Ass}}

\newcommand{\Spec}{\operatorname{Spec}}

\newcommand{\Ht}{\operatorname{ht}}

\newcommand{\Gfd}{\operatorname{Gfd}}
\newcommand{\V}{\operatorname{V}}

\newcommand{\Ext}{\operatorname{Ext}}
\newcommand{\Supp}{\operatorname{Supp}}

\newcommand{\Tor}{\operatorname{Tor}}
\newcommand{\Hom}{\operatorname{Hom}}

\newcommand{\Ann}{\operatorname{Ann}}
\newcommand{\Rad}{\operatorname{Rad}}

\newcommand{\Ndim}{\operatorname{Ndim}}

\newcommand{\depth}{\operatorname{depth}}
\newcommand{\width}{\operatorname{width}}
\newcommand{\Coass}{\operatorname{Coass}}

\newcommand{\Cosupp}{\operatorname{Cosupp}}

\newcommand{\adj}{\operatorname{adj}}

\newcommand{\vpl}{\operatornamewithlimits{\varprojlim}}

\newcommand{\lo}{\longrightarrow}
\newcommand{\fm}{\frak{m}}
\newcommand{\fp}{\frak{p}}

\newcommand{\fa}{\frak{a}}
\newcommand{\fb}{\frak{b}}

\newenvironment{prf}[1][Proof]
{\begin{proof}[\bf #1]}{\end{proof}}

\begin{document}

\author[M. Hatamkhani ]{Marziyeh Hatamkhani }
\title[Dual of Faltings' Theorems on Finiteness of Local Cohomology]{Dual of Faltings' Theorems on Finiteness of Local Cohomology}

\address{M. Hatamkhani, Department of Mathematics, Faculty of Science, Arak 
University, Arak, 38156-8-8349, Iran.} \email{m-hatamkhani@araku.ac.ir}



\begin{abstract}
 Let $R$ be a commutative Noetherian ring and $\fa$ an ideal of $R$.
We intend to establish the dual of two Faltings' Theorems for local homology modules of an Artinian module.
As a consequence of this, we show that, if $A$ is an Artinian module over semi-local complete ring $R$ and $j$ is an integer such that
$H_i^{\fa}(A)$ is Artinian for all $i<j$, then the set $\Coass_R(H_j^{\fa}(A))$ is finite.
We also introduce the notion of the $n$th Artinianness dimension $g_n^\fa(A)=\inf\{g^{\fa R_{\fp}}(^\fp A):  \fp\in\Spec(R) \ \  \text{and} \ \ \dim R/\fp\geq n\}$, for all $n\in\mathbb{N}_{0}$ and  prove that  $g_1^\fa(A)=\inf\{i\in\mathbb{N}_0: H_i^\fa(A)  \ \ \text{is not minimax}\}$, whenever $R$ is a semi-local complete ring. Moreover, in this situation we show that $\Coass_R(H_{g_1^\fa(A)}^\fa(A))$  is a finite set.

\end{abstract}

\maketitle

\section{Introduction}

Throughout this paper, $R$ denotes a commutative Noetherian ring with
nonzero identity, $\fa$ an ideal of $R$ and $M$  an $R$-module.

The theory of local cohomology was introduced by A. Grothendieck and developed so much after six decays
 by different commutative algebraists. While the local cohomology functors are studied in great detail not so
 much is known about the local homology functors.
 
The local homology functors were first studied by Matlis \cite{Mat} in 1974. Then the study
of this theory was continued by  work of Greenlees and May \cite{GM} and Alonso Tarr\'{i}o, Jerem\'{i}as L\'{o}pez
and Lipman \cite {AJL}. In \cite{HD}, the vanishing of local homology modules and also, in  \cite{CN} and in \cite{Ta}, the local homology of Artinian modules are studied.
It is shown that  local homology of Artinian modules admit a theory parallel to the theory for local cohomology
of finitely generated modules.

In section 2, we give some preliminaries on co-localization of an $R$- module and also local homology modules
which are necessary to the next sections.

 Let $M$ be a finitely generated $R$-module and let $\fa$ be an ideal of $R$.
An important theorem in  local cohomology is Faltings' Local- Global principle for the finiteness dimension of
 local cohomology modules \cite{F1}, which states that for a positive integer r, the $R_{\fp}$-module $H^i_{\fa R_{\fp}}(M_{\fp})$
 is finitely generated for all $i\leq r$ and for all $\fp\in Spec(R)$ if and only if the $R$-module $H^i_{\fa}(M)$ is 
 finitely generated for all $i\leq r$.
 
Recall that the finiteness dimension $f_{\fa}(M)$ of $M$ relative to $\fa$ is defined by  
$$f_{\fa}(M):=\inf\{i\in\mathbb{N}_0: H_{\fa}^i(M)\ \text{ is not finitely generated}\}=\inf\{i\in\mathbb{N}_0: \fa\nsubseteq \Rad({0:_RH_{\fa}^i(M)})\},$$
with the usual convention that the infimum of the empty set of integers is interpreted as $\infty$.

Let $\fb$ a second ideal of $R$ such that $\fb\subseteq\fa$. The $\fb$-finiteness dimension $f_{\fa}^{\fb}(M)$ of $M$ relative to $\fa$ is defined by $$f_{\fa}^{\fb}(M):=\inf\{i\in\mathbb{N}_0: \fb\nsubseteq \Rad({0:_RH_{\fa}^i(M)})\}.$$
Note that $f_{\fa}^{\fa}(M)=f_{\fa}(M)$. 

For a prime ideal $\fp\in\Spec(R)\setminus\V(\fa)$, the  $\fa$-adjusted depth of $M$ at $\fp$, denoted $\adj_{\fa}\depth M_{\fp}$, is defined by
$$\adj_{\fa}\depth M_{\fp}:=\depth M_{\fp}+\Ht(({\fa+\fp})/{\fp}).$$
 The  $\fb$-minimum  $\fa$-adjusted depth of $M$, denoted by $\lambda_{\fa}^{\fb}(M)$, is defined by 
$$\lambda_{\fa}^{\fb}(M):=\inf\{\adj_{\fa}\depth M_{\fp}: \fp\in\Spec(R)\setminus\V(\fb)\}$$$$=\inf\{\depth M_{\fp}+\Ht(({\fa+\fp})/{\fp}):\fp\in\Spec(R)\setminus\V(\fb)\}.$$

Another important theorem is Faltings' Annihilator Theorem, asserts that under milds restrictions on $R$, the invariants 
$\lambda_{\fa}^{\fb}(M)$ and  $f_{\fa}^{\fb}(M)$ are equal, see \cite{F2}.  

In section 3, we intoduce the concept  of the Artinianness dimension  $g^{\fa}(A)$ of  an Artinian $R$- module $A$ relative to $\fa$ by $$g^{\fa}(A):=\inf\{i\in\mathbb{N}: H_i^{\fa}(A)\ \text{ is not Artinian}\},$$ which is in some sence dual to the finiteness dimension and  prove that a dual of two Faltings' Theorems for local homology modules of an Artinian module over a semi-local complete ring.

As a consequence of this, we show that, if $A$ is an Artinian module over semi-local complete ring $R$ and $j$ is an integer such that
$H_i^{\fa}(A)$ is Artinian for all $i<j$, then the set $\Coass_R(H_j^{\fa}(A))$ is finite. 

 For a non-negative integer $n$, Bahmanpour et al., in \cite{BNS} introduced the notion of the $n$th finiteness dimension $f_{\fa}^n(M)$ of $M$ relative to $\fa$ by $$f_{\fa}^n(M):=\inf\{f_{\fa R_{\fp}}(M_{\fp}) : \fp\in\Supp(M/{\fa M}) \ \ \text{and} \ \ \dim R/\fp\geq n\}.$$
In the last section, we define the $n$th Artinianness dimension $g_n^\fa(A)$ of $A$ relative to $\fa$ by 
$$g_n^\fa(A)=\inf\{g^{\fa R_{\fp}}(^\fp A):  \fp\in\Spec(R) \ \  \text{and} \dim R/\fp\geq n\}$$
and show that the least integer $i$ such that  $H_i^{\fa}(A)$ is not minimax, equals to $g_1^\fa(A)$.
This implies that $\Coass_R(H_{g_1^\fa(A)}^\fa(A))$  is a finite set, whenever $R$ is semi-local complete ring.

For any ideal $\fb$ of $R$, the radical of $\fb$, denoted by $\Rad(\fb)$, is defined to be the set
$\{x\in R: x^n\in \fb \ \text{for some}\ \ n\in\mathbb{N}\}$. For any unexplained notation and terminology, we refer the reader to \cite{BS} and  \cite{GM}.

\section{Preliminaries}

Let $E_R$ be the module $\oplus E(R/\fm)$, where $\fm$ runs over the set of maximal ideals of $R$, and let 
$D_R(-)$ be the functor $\Hom_R(-,E_R)$.

\begin{definition}(See \cite{Ri})
Let $S$ be a multiplicatively closed subset of $R$. For any $R$-module $M$, the co-localization of $M$
relative to $S$ is the $S^{-1}R$-module $S_{-1}M=D_{S^{-1}R}(S^{-1}D_R(M))$.
If $S=R-\fp$ for some $\fp\in \Spec(R)$, we write $^{\fp}M$ for $S_{-1}M$.
\end{definition}
Note, as a composite of exact, additive functors, $S_{-1}(-)$ is exact and additive.

\begin{theorem}
Let $R$ be a semi-local and complete ring. If $S^{-1}R$ is also semi-local, but not necessarily complete,
then $S_{-1}(-)$ takes Artinian  $R$-modules to  Artinian $S^{-1}R$-modules.
\end{theorem}
\begin{prf}
See \cite [Theorem 2.3]{Ri}.
\end{prf}

\begin{lemma}
Let $R$  be semi-local with the Jacobson radical $J$ and $M$ an  $R$-module.
\begin{enumerate}
\item[i)] If $M$ is finitely generated  $R$-module, then $D_R(M)$  is an Artinian  $R$-module.
\item[ii)] Assume that $R$ is $J$- adically complete. If $M$ is an Artinian  $R$-module, then $D_R(M)$  is finitely generated  $R$-module.
\item[iii)] If $M$  is finitely generated  $R$-module, then $D_R(D_R(M))$  is isomorphic to $J$- adic completion of $M$. 
\item[iv)] If $M$ is an Artinian $R$-module, then  $D_R(D_R(M))$ is isomorphic to $M$.
\item[v)] For any  $R$-module $M$, $(0 :_R M)=(0 :_R D_R(M))$.
\end{enumerate}
\end{lemma}
\begin{prf}
See \cite[Theorem 1.6]{Oo}.
\end{prf}

The following lemma will be needed in Lemma 3.6.
\begin{lemma}
Let $R$ be semi-local complete and $S$ be a multiplicatively closed subset of $R$ such that the ring $S^{-1}R$ is semi-local. Then $$D_{S^{-1}R}(S_{-1}A)\cong S^{-1}(D_R(A))\otimes_{S^{-1}R}T,$$ where $T$ is Jacobson- adic completion of $S^{-1}R$.
\end{lemma} 
\begin{prf}
One has
 $$D_{S^{-1}R}(S_{-1}A)=\Hom_{S^{-1}R}(S_{-1}A,E_{S^{-1}R})=\Hom_{S^{-1}R}(\Hom_{S^{-1}R}(S^{-1}(\Hom_R(A,E_R)),E_{S^{-1}R}),E_{S^{-1}R}),$$
but $\Hom_R(A,E_R)$ is a finitely generated $R$-module and hence $S^{-1}(\Hom_R(A,E_R))$ is a finitely generated $S^{-1}R$-module.
Therefore, by Lemma 2.3(iii),  $D_{S^{-1}R}(S_{-1}A)\cong  S^{-1}(D_R(A))\otimes_{S^{-1}R}T$, which proves the result.
\end{prf}

Let $\fa$ be an ideal of $R$ and $\mathcal{C}_0(R)$ denote the
category of $R$-modules and $R$-homomorphisms. It is known that the
$\fa$-adic completion functor
$$\Lambda_{\fa}(-):=\underset{n}{\vpl}(R/\fa^n\otimes_R-):
\mathcal{C}_0(R)\lo \mathcal{C}_0(R)$$ is not right exact in
general. For any integer $i$, the $i$-th local homology functor with
respect to $\fa$ is defined as $i$-th left derived functor of
$\Lambda_{\fa}(-)$ and denoted by $H_i^{\fa}(-)$.

\begin{remark} Let $\phi: R\lo S$   be a homomorphism of commutative Noetherian rings.
For any ideal $\fa$ of $R$ and any  $S$-module $M$, we have $H_i^{\fa}(M)\cong H_i^{\fa S}(M)$  for all $i$.
\end{remark}
\begin{prf}
See \cite[Proposition 3.6]{Ri}.
\end{prf}
In what follows, we denote by $A$ an Artinian $R$-module.
\begin{theorem}
Let $R$ be semi-local. Then for any ideal $\fa$ of $R$, $H_i^{\fa}(A)\cong D_R(H^i_{\fa}(D_R(A)))$.
\end{theorem}
\begin{prf}
One has $H_i^{\fa}(A)\overset{(a)}{\cong}H_i^{\fa}(D_R(D_R(A)))\overset{(b)}{\cong}D_R(H^i_{\fa}(D_R(A)))$, 
where $(a)$  follows by Lemma 2.3(iv) and $(b)$ follows from \cite[Proposition 3.1]{Ri}.
\end{prf}

\section{Dual of Faltings' Theorems}
 Let $R$ be a commutative Noetherian semi-local complete ring. 
In this section we state and prove the dual of two Faltings' Theorem for local homology modules of an Artinian module.

The following Theorem is due to Coung, Nam\cite[Proposition 4.7]{CN}, without assumption semi-local completeness
of $R$. We give here a different proof by using Matlis duality.

\begin{theorem}
Let $R$ be semi-local complete. Then the following statements are equivalent for all $t\in\mathbb{N}$.
\begin{enumerate}
\item[i)] $H_i^{\fa}(A)$ is Artinian for all $i<t$;
\item[ii)] $\fa\subseteq \Rad({0:_RH_i^{\fa}(A)})$ for all $i<t$.
\end{enumerate}
\end{theorem}
\begin{prf}
By the above theorem it suffices to show that $D_R(H^i_{\fa}(D_R(A)))$ is Artinian for all $i<t$
 if and only if $\fa\subseteq \Rad({0:_RD_R(H^i_{\fa}(D_R(A)))})$. Lemma 2.3 now yeilds $D_R(H^i_{\fa}(D_R(A)))$  is Artinian 
for all $i<t$ if and only if $H^i_{\fa}(D_R(A))$ is finitely generated for all $i<t$ and by \cite[Proposition 9.12]{BS}, 
it is equivalent to $\fa\subseteq \Rad({0:_RH^i_{\fa}(D_R(A))})$ for all $i<t$.
But, from Lemma 2.3 we have $0:_RH^i_{\fa}(D_R(A))=0:_RD_R(H^i_{\fa}(D_R(A)))$, which completes the proof.
\end{prf}

\begin{definition}
Let $A$ be an Artinian $R$-module. In light of Theorem 3.1, we define the Artinianness dimension
 $g^{\fa}(A)$ of $A$ relative to $\fa$ by $$g^{\fa}(A):=\inf\{i\in\mathbb{N}_0: H_i^{\fa}(A)\ \text{ is not Artinian}\}=\inf\{i\in\mathbb{N}_0: \fa\nsubseteq \Rad({0:_RH_i^{\fa}(A)})\}.$$
\end{definition}

\begin{definition}
Let $\fb$ be a second ideal of $R$ such that $\fb\subseteq\fa$ .We define the $\fb$-artinianness dimension 
of $A$ relative to $\fa$ by $g^{\fa}_{\fb}(A):=\inf\{i\in\mathbb{N}_0: \fb\nsubseteq \Rad({0:_RH_i^{\fa}(A)})\}.$
Note that, $g^{\fa}_{\fa}(A)=g^{\fa}(A)$.
\end{definition}

\begin{lemma}
Let $R$ be semi-local complete. Then $g^{\fa}_{\fb}(A)=f_{\fa}^{\fb}(D_R(A))$.
\end{lemma}
\begin{prf}
It follows from Theorem 3.1 and Lemma 2.3(v).
\end{prf}

Let $M$ be an $R$-module. Recall that a prime ideal $\fp$ of $R$ is said to be a
coassociated prime ideal of $M$ if there is an Artinian quotient $L$
of $M$ such that $\fp=(0:_RL)$. The set of all coassociated prime
ideals of $M$ is denoted by $\Coass_RM$. (see \cite{Ya}). Also, $\mathcal{A}tt_RM$ is
defined by
$$\mathcal{A}tt_RM:=\{\fp\in \Spec(R)| \ \fp=(0:_RL) \  \ \text{for
some quotient} \   \ L \ \text{of} \  \  M \}.$$ Clearly,
$\Coass_RM\subseteq \mathcal{A}tt_RM$ and the equality holds if
either $R$ or $M$ is Artinian. More generally, if $M$ is representable, then 
it is easy to check that $\Coass_RM=\mathcal{A}tt_RM$.

\begin{theorem}
Let $R$ be semi-local complete. If $i\leq g^{\fa}(A)$, then $\Coass_R(H_i^{\fa}(A))$ is a finite set.
\end{theorem}
\begin{prf}
It follows from Lemma 3.4, \cite[Corollary 2.3]{BL} and \cite[Proposition 2.7]{Oo}.
\end{prf}

\begin{lemma}
Let the situation be as in the Definition 3.3, and let $S$ be a multiplicatively closed subset of complete semi-local ring $R$.
If $S^{-1}R$ is also semi-local, then $$g_{S^{-1}\fb}^{S^{-1}\fa}(S_{-1}A)=f_{S^{-1}\fa}^{S^{-1}\fb}(S^{-1}(D_R(A))).$$ 
\end{lemma}
\begin{prf}
We observe that $$(0:_{S^{-1}R}H_i^{S^{-1}\fa}(S_{-1}A))\overset{(a)}{=}(0:_{S^{-1}R}D_{S^{-1}R}(H^i_{S^{-1}\fa}(D_{S^{-1}R}(S_{-1}A))))
\overset{(b)}{=}(0:_{S^{-1}R}H^i_{S^{-1}\fa}(D_{S^{-1}R}(S_{-1}A)))$$$$\overset{(c)}{\cong} (0:_{S^{-1}R}H^i_{S^{-1}\fa}( S^{-1}(D_R(A)))\otimes_{S^{-1}R}T)=(0:_{S^{-1}R}H^i_{S^{-1}\fa}(S^{-1}(D_R(A)))),$$
where $T$ is Jacobson- adic completion of $S^{-1}R$ and $(a)$  follows by  Theorem 2.6, $(b)$ is by Lemma 2.3(v) and $(c)$ follows from Lemma 2.4 and Flat Base Change Theorem\cite[Theorem 4.3.2]{BS}. This yeilds the result.
\end{prf}

\begin{corollary}
Under the hypotheses of the above lemma, $g^{\fa}_{\fb}(A)\leq g_{S^{-1}\fb}^{S^{-1}\fa}(S_{-1}A)$. 
\end{corollary}
\begin{prf}
It follows by Lemma 3.4 , 3.6 and \cite[Exercise 9.16]{BS}. 
\end{prf}

\begin{theorem}(Local-global Principle for Artinianness Dimension)
Let $R$ be a semi-local complete ring. Then $g^{\fa}(A)=\inf \{g^{\fa R_{\fp}}(^{\fp}A) : \fp\in \Spec(R)\}.$
\end{theorem}
\begin{prf}
Note that $g^{\fa}(A)=f_{\fa}(D_R(A))$ and $g^{\fa R_{\fp}}(^{\fp}A)=f_{\fa R_{\fp}}(D_R(A)_{\fp})$ by Lemma 3.4 and 3.6,
respectively. Next, the result follows by Local-global Principle for Finiteness Dimension \cite[Theorem 9.6.2]{BS}.
\end{prf}

\begin{definition}
Let $R$ be a commuttive Noetherian ring and $M$ be an $R$-module.
An element $x$ of $R$ is called $M$-coregular if $M=xM$. A sequence $x_1,...,x_n$ in $R$ is called $M$-coregular sequence if 
\begin{enumerate}
\item[i)] $(0:_M(x_1,...,x_n))\neq 0$,
\item[ii)]  $x_i$ is an  $(0:_M(x_1,...,x_{i-1}))$-coregular element for each $i$, $i=1,...,n$.
\end{enumerate}
\end{definition}

\begin{definition}
Let $\fa$ be an ideal of $R$. The length of the longest  $M$-coregular sequence in $\fa$, is denoted by $\width_{\fa}(M)$.
If such sequence don't exist, then we write $\width_{\fa}(M)=\infty$.
\end{definition}

From \cite[Corollary 3.13]{Oo} and Lemma 2.3(iv), we have the following immediate consequence.
\begin{corollary}
Let $R$ be semi-local and let $A$ be an Artinian $R$-module. Then $\width_{\fa}(A)=\depth_{\fa}(D_R(A))$.
\end{corollary}

\begin{definition}
Let $A$ be an Artinian  $R$-module and $\fa$ an ideal of $R$. 
For a prime ideal $\fp\in\Spec(R)\setminus\V(\fa)$ the  $\fa$-adjusted width of $A$ at $\fp$, denoted $\adj_{\fa}\width^{\fp}A$, is defined by
$$\adj_{\fa}\width^{\fp}A:=\width^{\fp}A+\Ht(({\fa+\fp})/{\fp}).$$
Let $\fb$ be a second ideal of $R$ such that $\fb\subseteq\fa$ .We define the  $\fb$-minimum  $\fa$-adjusted width of $A$, denoted by $\gamma_{\fa}^{\fb}(A)$, by 
$$\gamma_{\fa}^{\fb}(A):=\inf\{\adj_{\fa}\width^{\fp}A: \fp\in\Spec(R)\setminus\V(\fb)\}$$$$=\inf\{\width^{\fp}A+\Ht(({\fa+\fp})/{\fp}):\fp\in\Spec(R)\setminus\V(\fb)\}$$
\end{definition}

\begin{corollary}
Let $R$ be semi-local complete and let $A$ be an Artinian $R$-module. Then $\gamma_{\fa}^{\fb}(A)=\lambda_{\fa}^{\fb}(D_R(A))$.
\end{corollary}
\begin{prf}
Follows by Corollary 3.11 and Lemma 2.4.
\end{prf}

We are now ready to state and prove the following main result.
\begin{theorem}
Let $R$ be a semi-local complete ring and let $A$ be an Artinian $R$-module. 
Suppose $\fb\subseteq\fa$ are two ideals of $R$. Then $g^{\fa}_{\fb}(A)=\gamma_{\fa}^{\fb}(A)$. 
\end{theorem}
\begin{prf}
First note that $D_R(A)$ is a finitely generated $R$-module. Now one has the equality 
$f_{\fa}^{\fb}(D_R(A))=\lambda_{\fa}^{\fb}(D_R(A))$ by \cite[Corollary 2.12]{KS}. 
Therefore the result follows by Corollary 3.13 and Lemma 3.4.
\end{prf}

\section{minimaxness and coartinianness of local homology modules}
Let $R$ be a commutative Noetherian ring, $\fa$ an ideal of $R$, and $M$ an $R$- module.
Recall that an $R$- module $M$ is called $\fa$- cofinite if $\Supp(M)\subseteq \V(\fa)$ and $\Ext^i_R(R/\fa,M)$ is 
finitely generated for all $i\geq 0$. The concept of $\fa$- cofinite modules were introduced by Hartshorne \cite{Ha} and has been proved
to be an important tool in the study of modules over a commutativt ring. 
Moreover, an $R$- module $M$ is said to be minimax, if there exists a finitely generated submodule $N$ of $M$, such that $M/N$
is Artinian. The class of minimax modules was introduced by Z\"{o}schinger \cite{Zo}.
 \\ 
 Bahmanpour et al., in \cite{BNS}  proved that $$f_{\fa}^1(M):=\inf\{i\in\mathbb{N}: H_{\fa}^i(M)\ \text{ is not minimax}\}.$$
\\
\par
In this section we  establish the dual of  the above equality for local homology modules.
The following concepts and lemmas are needed in the proof of these result.
\\In \cite{Ya}, Yassemi defined the co-support $\Cosupp_R(M)$ of an $R$- module $M$ to be the set of primes $\fp$ such that there exists
a cocyclic homomorphic image $L$ of $M$ with $\Ann(L)\subseteq \fp$. Recall that a module is cocyclic if it is a submodule of $E(\frac{R}{\fm})$ for some maximal ideal $\fm$ of $R$. It is clear that for any $R$- module $M$, $\Coass_(M)\subseteq\Cosupp_R(M)$.
\begin{lemma}(\cite{CN}[3.8] and \cite{Na}[3.8]).
Let $A$ be an Artinian $R$- module. Then $$\Cosupp_R(H_i^\fa(A))\subseteq\Cosupp_R(A)\bigcap \V(\fa)$$ for all $i\geq 0$.
\end{lemma}
\par
Nam introduced the following definition of coartinian modules which is in some sense dual to the conctpt of cofinite modules.

\begin{definition}(\cite{Na}[4.1])
Let $\fa$ be an ideal of $R$ and $M$ an $R$- module. $M$ is said to be $\fa$- coartinian if $\Cosupp_R(M)\subseteq \V(\fa)$ and $\Tor_i^R(R/\fa,M)$ is an artinian $R$- module for each $i$.
\end{definition}
\begin{remark}
In virtue of \cite{DFT}[Corollary 4.4], an $R$- module $M$ with $\Cosupp_R(M)\subseteq \V(\fa)$  is  $\fa$- coartinian if and only if
$\Ext^i_R(R/\fa,M)$ is an artinian $R$- module for each $i$.
\end{remark}
The following are some basic properties of coartinian modules.
\begin{lemma}(\cite{Na}[4.2(i)]
If $0\longrightarrow M'\longrightarrow M\longrightarrow M''\longrightarrow 0$ is an exact sequence and two of the modules in the 
exact sequence are $\fa$- coartinian, then so is the third one.
\end{lemma}

\begin{lemma}
Let $R$  be semi-local with the Jacobson radical $J$ and $M$ an  $R$-module.
\begin{enumerate}
\item[i)] If $D(M)$ is $\fa$- cofinite, then $M$ is $\fa$- coartinian.
\par
Moreover, if $R$ is $J$- adically complete,
\item[ii)]   If $M$ is $\fa$- coartinian, then $D(M)$ is $\fa$- cofinite.
\item[iii)] If $D(M)$ is $\fa$- coartinian, then $M$ is $\fa$- cofinite. 
\end{enumerate}
\end{lemma}
\begin{prf}
(i), (ii). By \cite{Ya}[Theorem 3.2], we have $\Supp(D(M))=\Cosupp_R(M)$, and by \cite{Ri}[Proposition 1.4],
 $$D(\Tor_i^R(R/\fa,M))\cong\Ext^i_R(R/\fa, D(M)).$$ Now the result follows from Lemma  (i),(ii).
 \par
(ii) It follows from \cite{Ya}[Theorem 3.2] and \cite{Ri}[Proposition 1.4] that $\Supp_R(M)\subseteq\Cosupp_R(D(M))$ and
$$\Tor_i^R(R/\fa,D(M))\cong D(\Ext^i_R(R/\fa, M)).$$
In virtue of Lemma 2.3, $\Ext^i_R(R/\fa, M)$ is finite if and only if $\Tor_i^R(R/\fa,D(M))$ is Artinian.
The proof is complete.
\end{prf}

\begin{theorem}
Let $R$  be semi-local with the Jacobson radical $J$ and $M$ an  $R$-module with $\Cosupp_R(M)\subseteq \V(\fa)$. 
Assume that $R$ is $J$- adically complete. Then $M$ is finitely generated and $\fa$- coartinian if and only if $M/\fa M$
has finite length.
\end{theorem}
\begin{prf}
 $(\Rightarrow)$ is clear.
 \par 
 $(\Leftarrow)$ By Lemma 2.3,$D_R(M/{\fa M})$ has finite length, but $D_R(M/{\fa M})\cong (0:_{D_R(M)}\fa)$
 and so  \cite{Me}[proposition 4.1] follows that $D_R(M)$ is Artinian and  $\fa$- cofinite. Now the result follows from Lemma 2.3 and 
4.5.
\end{prf}
\begin{theorem}
Let  $R$ be as Theorem 4.6 and $M$ be a minimax $R$- module with $\Cosupp_R(M)\subseteq \V(\fa)$. Then $M$ is $\fa$- coartinian if and only if $M/{\fa M}$
is Artinian.
\end{theorem} 
\begin{prf}
There exists a finitely generated submodule $N$ of $M$ such that $M/N$ is Artinian. The exact sequence 
$$0\longrightarrow N\longrightarrow M\longrightarrow M/N\longrightarrow 0$$ induces a long exact sequence 
$$\Tor_1^R(R/\fa,M/N)\longrightarrow R/\fa\otimes_RN\longrightarrow R/\fa\otimes_RM.$$
Since $M/N$ is  $\fa$- coartinian, we deduce from the above exact sequence that $R$-module $N/{\fa N}\cong R/\fa\otimes_RN$ is Artinian and so by assumption, it has finite length. We also have $\Cosupp_R(N)\subseteq\Cosupp_R(M)\subseteq \V(\fa)$.
Hence by Theorem 4.6, $N$ is  $\fa$- coartinian. Now the result follows from Lemma 4.4.
\end{prf}
\begin{definition}
Let $R$ be a commutative Noetherian ring and $\fa$ an ideal of $R$. Let $A$ be an Artinian $R$- module and $n$ be a non-negative
integer. We define the $n$th Artinianness dimension $g_n^\fa(A)$ of $A$ relative to $\fa$ by 
$$g_n^\fa(A)=\inf\{g^{\fa R_{\fp}}(^\fp A):  \fp\in\Spec(R) \ \  \text{and} \ \ \dim R/\fp\geq n\}.$$
Note that $g_n^{\fa}(A)$ is either a positive integer or $\infty$ and that $g_0^{\fa}(A)=g^{\fa}(A)$.
\end{definition}
\begin{lemma}
Let $R$ be a semilocal and $J$- adically complete ring and $A$ an Artinian $R$- module.
Then  $g_n^\fa(A)=f^n_\fa(D_R(A))$.
\end{lemma}
\begin{prf}
It follows from Lemma 3.4 and 2.4.
\end{prf}
\begin{corollary}
Let $R$ and $A$ be as in Lemma 4.9. Then $g_1^\fa(A)=\inf\{i\in\mathbb{N}_0: H_i^\fa(A)  \ \ \text{is not minimax}\}$.
\end{corollary}
\begin{prf}
Note that an $R$- module $M$ is minimax if and only if $D_R(M)$ is minimax. Now the result follows from Lemma 4.9, Theorem 2.6 and
\cite {BNS}[Corollary 2.4].
\end{prf}
\begin{theorem}
Let $R$ and $A$ be as in Lemma 4.9. Let $n$ be a non-negative integer such that $\Tor_j^R(R/\fa, H_i^\fa(A))$
 is Artinian for all $j\in\mathbb{Z}$ and $i<n$. Then the $R$- modules 
 $$\Tor_0^R(R/\fa, H_n^\fa(A))  \ \ \text{and} \ \ \Tor_1^R(R/\fa, H_n^\fa(A))$$ are Artinian.
\end{theorem}
\begin{prf}  
From Lemma 4.1 and assumption  $H_i^\fa(A)$ is $\fa$- coartinian for all  $i<n$ and so by 
Theorem 2.6  and Lemma 4.5(iii) $H^i_\fa(D_R(A))$ is $\fa$- cofinite for all $i<n$. It follows from \cite{Kh}[Corollary 3.5] that $R$- modules
$$\Hom_R(R/\fa, H^n_\fa(D_R(A))) \ \ \text{and} \ \  \Ext^1(R/\fa, H^n_\fa(D_R(A)))$$ are finitely generated.
Hence by Lemma 2.3 and  \cite{Ri}[Propositon 4.1] $\Tor_i^R(R/\fa, D_R(H^n_\fa(D_R(A))))$ is Artinian for $i=0, 1$.
Now the result follows from Theorem 2.6. 
\end{prf}
\begin{theorem}
Let $R$ and $A$ be as in Lemma 4.9. Then the following conditions holds:
\begin{enumerate}
\item[i)] The $R$- module $H_i^\fa(A)$ is minimax and $\fa$- coartinian for all $i<g_1^\fa(A)$;
\item[ii)] For all submodule $L$ of $H_{g_\fa^1(A)}^\fa(A)$ such that $H_{g_1^\fa(A)}^\fa(A)/L$ is minimax, the $R$- modules 
$$L/{\fa L} \ \ \text{and} \ \ \Tor_1^R(R/\fa, L)$$ are Artinian, whenever $g_1^\fa(A)$ is finite.
\end{enumerate}
\end{theorem}
\begin{prf}
(i) Let $t:=g_1^\fa(A)$, by Corollary 4.10 ,$H_i^\fa(A)$ is minimax for all $i<t$. Also by Lemma 4.9, we have $t=f^1_\fa(D_R(A))$ which follows from \cite {BNS}[Theorem 2.3] that $H^i_\fa(D_R(A))$ is $\fa$- cofinite and minimax for all $i<t$. 
Hence $H^i_\fa(D_R(A))$ is Matlis reflexive for all $i<t$ by \cite{Zo}. Therfore by Theorem 2.6,
we have $$D_R(H_i^\fa(A))\cong D_R(D_R(H^i_\fa(D_R(A))))\cong H^i_\fa(D_R(A)).$$
Now the result follows from Lemma 4.5(i).
\par
(ii) By Theorem 4.11 $$H_{g_1^\fa(A)}^\fa(A)/\fa H_{g_1^\fa(A)}^\fa(A) \ \ \text{and} \ \ \Tor_1^R(R/\fa, H_{g_1^\fa(A)}^\fa(A))$$
are Artinian. Let $N:=H_{g_1^\fa(A)}^\fa(A)/L$.  The exact sequence
 $$0\longrightarrow L\longrightarrow H_{g_1^\fa(A)}^\fa(A)\longrightarrow N\longrightarrow 0$$
 induces a long exact sequence $$\Tor_2^R(R/\fa, N)\longrightarrow \Tor_1^R(R/\fa, L)\longrightarrow \Tor_1^R(R/\fa, H_{g_1^\fa(A)}^\fa(A))\longrightarrow  \Tor_1^R(R/\fa, N)$$$$\longrightarrow R/\fa \otimes_R L\longrightarrow R/\fa \otimes_RH_{g_1^\fa(A)}^\fa(A)\longrightarrow R/\fa\otimes_R N\longrightarrow 0,$$
 which follows that $R/\fa\otimes_R N$ is Artinian. Since $N$ is minimax and 
 $\Cosupp_R(N)\subseteq \Cosupp_RH_{g_1^\fa(A)}^\fa(A)\subseteq \V(\fa)$, by Theorem 4.7, we have $N$ is $\fa$- coartinian.
 Thus the result follows from the above long exact sequence. 
\end{prf}
\begin{corollary}
Let $R$ and $A$ be as in Lemma 4.9. Then $\Coass_R(H_{g_1^\fa(A)}^\fa(A))$  is a finite set.
\end{corollary}
\begin{prf}
It follows from  Theorem 4.12, Lemma  4.1 and \cite{Ya}[Lemma 1.22].
\end{prf}



\end{document}